\documentclass[11pt, oneside]{article}   	
\usepackage{geometry}                		
\usepackage{graphicx}				
\usepackage{latexsym}
\usepackage{amsmath}
\usepackage{amsthm}
\usepackage{amscd}
\usepackage{amssymb}
\usepackage{mathrsfs}
\usepackage{mathabx}

\theoremstyle{definition}

\theoremstyle{remark}

\newtheorem*{Kazdan-Warner}{The Kazdan-Warner identity}

\newcommand{\tr}{\mathrm{tr}}

\renewcommand{\d}{\mathrm{d}}
\newcommand{\dvol}{\mathrm{dvol}}
\newcommand{\slashnabla}{\nabla\llap{/\,}}
\newcommand{\slashd}{\d\llap{/}}
\newcommand{\slashepsilon}{\epsilon}
\newcommand{\slashdiv}{\mathrm{d}\rlap\slash \mathrm{i}\mathrm{v}}
\newcommand{\slashcurl}{\mathrm{c}\rlap{\;/} \mathrm{ur} \mathrm{l}}
\newcommand{\slashDelta}{\Delta\llap{/\,}}
\newcommand{\sym}{\mathrm{sym}\,}
\numberwithin{equation}{section}

\title{A proof of the Kazdan-Warner identity via the Minkowski spacetime}
\author{Pengyu Le}
\date{}

\begin{document}
\maketitle

\begin{abstract}
Any $2$-dim Riemannian manifold with spherical topology can be embedded isometrically into a lightcone of the Minkowski spacetime. We apply this fact to give a proof of the Kazdan-Warner identity.
\end{abstract}

\section{The Kazdan-Warner identity}\label{The Kazdan-Warner identity}

$(\mathbb{S}^2, g_{\mathbb{S}^2})$ is the sphere of radius $1$ at the origin in the $3$-dim Euclidean space $\mathbb{E}^3$. $(\mathbb{S}^2,g_{\mathbb{S}^2})$ has the constant Gauss curvature $1$. $\{x_1,x_2,x_3\}$ is the rectangular coordinate system of $\mathbb{E}^3$. We denote $\nabla$ the Levi-Civita connection on $(\mathbb{S}^2,g_{\mathbb{S}^2})$. Then the vector fields $\nabla x_1, \nabla x_2, \nabla x_3$ are conformal Killing vector fields on $(\mathbb{S}^2,g_{\mathbb{S}^2})$, i.e. the diffeomorphisms generated by them are conformal.

$g$ is another Riemannian metric on the sphere $\mathbb{S}^2$. By the uniformization theorem, there exists a function $f$ on $\mathbb{S}^2$ such that  the conformal metric $e^{-2f}g$ has the constant Gauss curvature 1. Hence there is a diffeomorphism $\psi: \mathbb{S}^2 \rightarrow \mathbb{S}^2$ such that $\psi^{*}g=e^{2f\circ\psi} g_{\mathbb{S}^2}$. 
 
In the following of this note, we assume that $g$ is conformal to the standard metric $g_{\mathbb{S}^2}$ with the conformal factor $e^{2f}$. Let $K_g$ be the Gauss curvature of $(\mathbb{S}^2,g)$. Then we have the following from \cite{K-W}:
\begin{Kazdan-Warner}
\begin{equation}
\int_{\mathbb{S}^2} \langle \nabla K_g, \nabla x_i \rangle_{g_{\mathbb{S}^2}} \dvol_g =0.
\end{equation}
\end{Kazdan-Warner}

We denote $\tilde{\nabla}$ the Levi-Civita connection on $(\mathbb{S}^2,g)$. Then we can rewrite the Kazdan-Warner identity as following
\begin{equation}
\int_{\mathbb{S}^2} \langle \tilde\nabla K_g, \nabla x_i \rangle_g \dvol_g =0.
\end{equation}
Since $\nabla x_i$ is a conformal Killing vector field for $g_{\mathbb{S}^2}$, it is also a conformal Killing vector field for $g$. Actually we have that for any conformal Killing vector field $X$ on $(\mathbb{S}^2,g)$,
\begin{equation}\label{The generalized K-W identity}
\int_{\mathbb{S}^2} \langle \tilde{\nabla} K_g, X \rangle_g \dvol_g =0.
\end{equation}
We will give a proof of \eqref{The generalized K-W identity}, thus the Kazdan-Warner identity follows.

\section{The Minkowski spacetime}

$(\mathbb{M}, \eta)$ is the $4$-dim Minkowski spacetime. $\{x_0, x_1,x_2,x_3\}$ is the rectangular coordinate system of $\mathbb{M}$. The metric $\eta=-\d x_0^2 + \d x_1^2 + \d x_2^2 + \d x_3^2$.

$t,r$ are two functions on $\mathbb{M}$: $t=x_0$ and $r= \sqrt{x_1^2+x_2^2+x_3^2}$. We define the optical functions $u,v$ by $u=\frac{1}{2}(t-r)$ and $v=\frac{1}{2}(t+r)$. So $t=u+v$ and $r=v-u$.

Let $C_u$ be the level set of $u$, which is the future lightcone at the point $(u,0,0,0)$. Let $\underline{C}_v$ be the level set of $v$, which is the past lightcone at the point $(v,0,0,0)$. Let $S_{u,v}$ be the intersection of $C_u$ and $\underline{C}_v$. $S_{u,v}$ is a sphere of radius $r=v-u$.

We define a map $\phi : \mathbb{M}\setminus \{r= 0\} \rightarrow \mathbb{S}^2$ by $\phi: x=(x_0,x_1,x_2,x_3) \mapsto (\frac{x_1}{r},\frac{x_2}{r},\frac{x_3}{r})$. $\mathbb{R}^2_{v>u}$ is the open half-plane $\{(v,u)\in\mathbb{R}^2 \mid v>u\}$. Then we have another coordinate system $\Phi$ on $\mathbb{M}\setminus \{ r=0 \}$ given by
\begin{eqnarray}
\Phi :&   \mathbb{M}\setminus \{r=0\} &\rightarrow  \mathbb{R}^2_{v>u} \times \mathbb{S}^2, \nonumber
\\
&
x=(x_0,x_1,x_2,x_3)& \mapsto  \left( (v,u),\phi(x)=(\frac{x_1}{r},\frac{x_2}{r},\frac{x_3}{r})\right).
\end{eqnarray}

Let $\partial_t$ be the vector field $\partial_0$ and $\partial_r$ be the vector field $\frac{x_1}{r} \partial_1+\frac{x_2}{r} \partial_2+ \frac{x_3}{r} \partial_3$. In the coordinate system $\Phi$, the coordinate vector fields $\partial_v = \partial_t + \partial_r$ and $\partial_u = \partial_t - \partial_r$.

$\partial_u$ and $\partial_v$ are both null vector fields, i.e. $\eta(\partial_u,\partial_u)=\eta(\partial_v,\partial_v)=0$. The inner product of $\partial_u$ and $\partial_v$ is $-2$, i.e. $\eta(\partial_u,\partial_v)=-2$. We see that $\partial_u$ and $\partial_v$ are orthonormal to the tangent space of any $S_{u,v}$. So in the coordinate system $\Phi$, the metric $\eta=-2(\d u \otimes \d v + \d v \otimes \d u) + r^2 g_{\mathbb{S}^2}$.

In particular, when we restrict the coordinate system $\Phi$ on the lightcone $C_{0}\setminus \{o\}$, we get a coordinate system of $C_{0}\setminus \{o\}$: 
\begin{equation}
\Phi|_{C_0\setminus \{o\} }: \quad C_0\setminus \{o\} \rightarrow \mathbb{R}_{+} \times \mathbb{S}^2, \quad x\in C_0\setminus \{o\} \mapsto (v,\phi(x)).
\end{equation}
The induced metric $\eta|_{C_0\setminus \{o\} }=v^2 g_{\mathbb{S}^2}$ is degenerated.

\section{The isometric embedding of $(\mathbb{S}^2, g)$ into a lightcone of the Minkowski spacetime}\label{embedding}

Via the coordinate system $\Phi|_{C_0\setminus\{o\}}$, we can represent any closed spacelike surface in $C_0\setminus \{o\}$ as a graph of a function on $\mathbb{S}^2$. $S$ is a closed spacelike surface in $C_0\setminus\{o\}$, there exists a unique function $h$ on $\mathbb{S}^2$ such that
\begin{equation}
\Phi|_{C_0\setminus\{o\}}(S)= \left\{ (e^{h(\theta)}, \theta)\in \mathbb{R}_{+}\times\mathbb{S}^2 \mid \theta\in \mathbb{S}^2 \right\}.
\end{equation}

Conversely for any function $h$ on $\mathbb{S}^2$, the map
\begin{equation}\label{psi_h}
\psi_h: \quad \mathbb{S}^2 \rightarrow \mathbb{R}^+ \times \mathbb{S}^2, \quad \theta\in \mathbb{S}^2 \mapsto (e^{h(\theta)},\theta)
\end{equation}
is an embedding of $\mathbb{S}^2$ into $C_0\setminus \{o\}$. Its image
\begin{equation}\label{S_h}
S_h = \psi_h(\mathbb{S}^2) = \left\{ (e^{h(\theta)}, \theta) \in \mathbb{R}_{+}\times\mathbb{S}^2 \mid \theta\in \mathbb{S}^2 \right\}
\end{equation}
is a closed spacelike surface in $C_0\setminus \{o\}$. Moreover, $\psi_h^* \left( \eta|_{S_h} \right)= e^{2h} g_{\mathbb{S}^2}$ since $\eta|_{C_0\setminus\{o\}} =v^2 g_{\mathbb{S}^2}$.

Hence, we can embed $(\mathbb{S}^2,g=e^{2f} g_{\mathbb{S}^2})$ isometrically into $C_0\setminus \{o\}$ as above by taking $h=f$.

\section{The geometry of a spacelike surface in the Minkowski spacetime}\label{geo of spacelike surface}

Let $S$ be a orientable spacelike surface in the Minkowski spacetime $(\mathbb{M},\eta)$. $TS$ is the tangent bundle and $NS$ is the normal bundle of $S$ in $(\mathbb{M},\eta)$.

Since $S$ is spacelike, $\eta|_{TS}$ is positive definite and $\eta|_{NS}$ is of signature $(1,1)$. Hence we can choose a null frame $\{L,\underline{L}\}$ of $NS$ such that $L$ and $\underline{L}$ are both future-directed null vector fields and their inner product is $-2$, i.e.
\begin{equation}
\eta(L,L)=\eta(\underline{L},\underline{L})=0, \quad \eta(L,\underline{L})= -2.
\end{equation}
Such a choice isn't unique, since for any positive function $a$ on $S$, the frame $\{aL,a^{-1}\underline{L}\}$ also satisfies the above conditions.

In the following, we fix such a null frame $\{L,\underline{L}\}$. Then we can choose an orthonormal frame $\{e_1,e_2\}$ of $TS$ at least locally such that $\{e_1,e_2,e_3=\underline{L},e_4=L\}$ is positive oriented in $\mathbb{M}$. Choose the orientation of $S$ to be the orientation of $\{e_1,e_2\}$. Let $A,B=1,2$. 
 
The intrinsic geometry of $S$ is given by the induced metric $\eta|_{S}$. Let $\slashnabla$ be the Levi-Civita connection on $(S,\eta|_S)$ and $\slashd$ be the exterior derivative on $S$. Let $\slashepsilon$ be the volume form on $(S,\eta|_S)$. We define the intrinsic differential operators $\slashcurl,\slashdiv$ on $S$ by
\begin{eqnarray}
& \slashcurl\, \omega = \slashepsilon^{AB}\,\slashnabla_A \omega_B, &\text{ for any 1-form } \omega \text{ on } S;
\\
& (\slashdiv\, T)_A = \slashnabla^B T_{BA},  &\text{ for any symmetric 2-tensor field } T \text{ on } S.
\end{eqnarray}

The extrinsic geometry of $S$ in $(\mathbb{M},\eta)$ is given by the null second fundamental forms $\chi,\underline{\chi}$ and the torsion $\zeta$ defined as following:
\begin{eqnarray}
\chi(X,Y)=\eta(\nabla_X L,Y), & \qquad \text{for any } X,Y \in TS;
\\
\underline{\chi} (X,Y)=\eta(\nabla_X \underline{L},Y), & \qquad \text{for any } X,Y\in TS;
\\
\zeta(X) = \frac{1}{2} \eta(\nabla_X L, \underline{L}), & \qquad \text{for any } X \in TS. \label{torsion}
\end{eqnarray}
$\chi,\underline{\chi}$ are covariant symmetric 2-tensor fields and $\zeta$ is a 1-form on $S$. Let $\tr \chi$ and $\tr \underline{\chi}$ be the traces of $\chi$ and $\underline{\chi}$ on $S$, i.e. the contractions with $\eta|_S^{-1}$ on $TS$. Let $\hat{\chi}$ and $\hat{\underline{\chi}}$ be the tracefree parts of $\chi$ and $\underline{\chi}$ on $S$.

In analogy with the Gauss equation and Codazzi equation for a surface in $3$-dim Euclidean spacetime, we have the following equations:
\begin{description}
\item[The Gauss equations]
\begin{eqnarray}
K+\frac{1}{4}\tr\chi \, \tr\underline{\chi}-\frac{1}{2}(\hat{\chi},\hat{\underline{\chi}})_{\eta|S}=0,
\\
\slashcurl\, \zeta + \frac{1}{2} \hat{\underline{\chi}} \wedge\hat{\chi} =0;
\end{eqnarray}
\item[The Codazzi equations]
\begin{eqnarray}
\slashdiv\, \hat{\chi} -\frac{1}{2}\slashd\, \tr \chi + \hat{\chi} \cdot \zeta -\frac{1}{2} \tr\chi\, \zeta=0,
\\
\slashdiv\, \hat{\underline{\chi}} -\frac{1}{2} \slashd\, \tr\underline{\chi} -\hat{\underline{\chi}}\cdot \zeta +\frac{1}{2} \tr\underline{\chi}\, \zeta=0;
\end{eqnarray}
\end{description}
where
\begin{equation}
(\hat{\chi},\hat{\underline{\chi}})_{\eta|_S} = (\eta|_S^{-1})^{AC}(\eta|_S^{-1})^{BD} \hat{\chi}_{AB}\hat{\underline{\chi}}_{CD}, 
\quad
\hat{\underline{\chi}}\wedge \hat{\chi} =\slashepsilon^{AB} \hat{\underline{\chi}}_A^{\phantom{A} C} \hat{\chi}_{CB},
\end{equation}
and
\begin{equation}
(\hat{\chi} \cdot \zeta)_A =\hat{\chi}_A^{\phantom{A}B} \zeta_B, \quad (\hat{\underline{\chi}}\cdot \zeta)_A= \hat{\underline{\chi}}_A^{\phantom{A}B}\zeta_B.
\end{equation}
One can find the proofs of these equations in \cite{C}.

We can apply these equations to $S$ in the lightcone $C_0\setminus \{o\}$. We see $\partial_v$ along $S$ is a null vector field in $NS$. So we choose $L$ to be the null vector field $\partial_v$ over $S$. Then we can find $\underline{L}$ in $NS$ such that $\{L,\underline{L}\}$ is a null frame of $NS$.

Since the metric $\eta|_{C_0\setminus \{o\}} = v^2 g_{\mathbb{S}^2}$, we see that the induced metric on $S$ deforms conformally when we deform $S$ in the direction of $\partial_v$. This means that the seconded fundamental form $\chi$ of $S$ is a multiple of the induced metric $\eta|_{S}$. Hence $\hat{\chi}=0$. So on $S$, the equations are simpler:
\begin{eqnarray}
K+\frac{1}{4}\tr\chi \tr\underline{\chi}=0,  \label{Gauss 1}
\\
 \slashcurl \, \zeta =0, \label{Gauss 2}
\\
\slashd \,\tr\chi + \tr\chi \, \zeta =0, \label{Codazzi 1}
\\
\slashdiv \, \hat{\underline{\chi}} -\frac{1}{2} \slashd\, \tr\underline{\chi} -\hat{\underline{\chi}} \cdot \zeta +\frac{1}{2} \tr\underline{\chi} \,  \zeta=0.  \label{Codazzi 2}
\end{eqnarray}

\section{The proof of the Kazdan-Warner identity}

$S$ is a closed surface in $C_0\setminus \{o\}$. Let $X$ be a conformal Killling vector field on $(S,\eta|_S)$. The deformation tensor field $^{(X)}\pi=\mathcal{L}_X \eta|_S$ of $X$ is a multiple of $\eta|_S$. Then $^{(X)}\pi=^{(X)}\Omega \cdot \eta|_S$ where $^{(X)}\Omega=\frac{1}{2}\slashdiv X$. We define the operator $\sym$ by
\begin{equation}
(\sym T)_{AB} = T_{AB}+T_{BA}, \text{ for any covariant 2-tensor field } T \text{ on } S. 
\end{equation}
Then $\sym (\slashnabla X) =\mathcal{L}_X \eta|_S = ^{(X)}\pi =\Omega^X \cdot \eta|_S$.

\begin{eqnarray*}
\int_{S} \langle \slashnabla K, X \rangle_{\eta|_{S}} \dvol_{\eta|_{S}}
&\stackrel{\eqref{Gauss 1}}{=}&
\int_{S} - \frac{1}{4} \langle \slashnabla(\tr\chi \, \tr\underline\chi), X \rangle_{\eta|_{S}} \dvol_{\eta|_{S}}
\\
&=&
\int_{S} -\frac{1}{4}  \langle \tr\underline\chi \slashnabla \tr\chi +  \tr\chi \slashnabla \tr\underline\chi, X \rangle_{\eta|_{S}} \dvol_{\eta|_{S}} 
\\
&\stackrel{\eqref{Codazzi 1}\eqref{Codazzi 2}}{=}&
\int_{S} -\frac{1}{2} \left\{ \tr\chi \, \slashdiv\, \hat{\underline\chi} \cdot X - \hat{\underline{\chi}} (\tr\chi \, \zeta, X ) \right\} \dvol_{\eta|_{S} }
\\
&=&
 \int_{S} \frac{1}{2}\left\{  \tr\chi\, \langle \hat{\underline{\chi}}, \slashnabla X\rangle_{\eta|_{S}}  + \hat{\underline{\chi}}(\slashnabla \tr\chi+ \tr\chi \, \zeta, X)  \right\} \dvol_{\eta|_{S}} 
\\
&\stackrel{\eqref{Codazzi 1}}{=}&
\int_{S} \frac{1}{4} \tr\chi\, \langle \hat{\underline{\chi}} , \sym (\slashnabla X)\rangle_{\eta|_{S}} \dvol_{\eta|_{S}}
\\
&=&
\int_{S} \frac{1}{4} \tr\chi\, \langle \hat{\underline{\chi}} , \Omega_X \cdot \eta|_S\rangle_{\eta|_{S}} \dvol_{\eta|_{S}}
\\
&=&
0.
\end{eqnarray*}
The last equality follows from that $\hat{\underline{\chi}}$ is tracefree.

Together with the constructions in section~\ref{embedding}, we prove \eqref{The generalized K-W identity}.

\section{The gauge transformations on the normal bundle of a spacelike surface in the Minkowski spacetime}

Recall that in the section~\ref{geo of spacelike surface}, we introduced the normal bundle $NS$ of a oriented spacelike surface $S$ in the Minkowski spacetime $(\mathbb{M},\eta)$. Since that for any $p\in S$, the normal space $N_p S$ endowed with the induced metric $\eta|_{N_p S}$ is isometric to the $2$-dim Minkowski spacetime, we have a $1$-dim non-compact abelian group of isometries for $(N_p S,\eta|_{N_p S})$. The group of isometries on $(N_p S, \eta|_{N_p S})$ is just the group of Lorentz rotations of the $2$-dim Minkowski spacetime. We can explicitly write down the isometries via the null frame $\{L,\underline{L}\}$. Any positive number $a\in \mathbb{R}_{>0}$, we have the mapping $\mathcal{L}_a: N_p S \rightarrow N_p S$ defined by
\begin{equation}
\mathcal{L}_a: \quad L_p \rightarrow aL_p, \quad \underline{L}_p \rightarrow a^{-1}\underline{L}_p.
\end{equation}
The group structure is simply given by $\mathcal{L}_a \circ \mathcal{L}_b = \mathcal{L}_{ab}$ for any $a,b\in \mathbb{R}_{>0}$.

So the normal bundle $(NS,\eta|_{NS})$ is a vector bundle on $S$ with the group action of $(\mathbb{R}_{>0},\cdot)$. Moreover, the null frame bundle of $(NS,\eta|_{NS})$ is a principal $(\mathbb{R}_{>0},\cdot)$-bundle. This principal bundle is actually trivial since we can find a global section, which is a global null frame. The parallel transport on the normal bundle $NS$ defines a principal connection on the null frame bundle. We see that the torsion $\zeta$ of a null frame $\{L,\underline{L}\}$ is actually the connection 1-form for this principal connection. Assume now that $a$ is a positive function over $S$, then $\{ aL,a^{-1}\underline{L}\}$ is another null frame of $(NS,\eta|_{NS})$. Let us denote $\zeta_a$ being the torsion for $\{aL,a^{-1}\underline{L}\}$. Direct calculation by the definition of the torsion \eqref{torsion} shows that
\begin{equation}
\zeta_a = \zeta - a^{-1} \slashd \, a = \zeta -\slashd \, \log a,
\end{equation}
which is just the transformation formula for the connection form. However $\slashcurl \zeta_a =\slashcurl \zeta$ keeps invariant, because it is actually the curvature of this connection. In particular, we can choose a positive function $a$ such that
\begin{equation}
\slashdiv \, \zeta_a = \slashdiv \,\zeta - \slashDelta \log a =0.
\end{equation}

Now Assume that the oriented spacelike surface $S$ is contained in the lightcone $C_0\setminus \{o\}$. We take the null frame $\{L,\underline{L}\}$ of $(NS_h,\eta|_{NS_h})$ such that its torsion $\zeta$ satisfies $\slashdiv \, \zeta = 0$. Associated with this null frame $\{L,\underline{L}\}$, we have the Gauss equations and Codazzi equations. Since $\hat{\chi}=0$ still holds, we have $\slashcurl \, \zeta =0$. Then the equations
\begin{equation}
\slashdiv \, \zeta =0, \quad \slashcurl \, \zeta =0,
\end{equation}
imply that $\zeta =0$. Then we have
\begin{eqnarray}
K+\frac{1}{4} \tr\chi \tr\underline\chi =0,
\\
\slashd \, \tr \chi =0,
\\
\slashdiv \, \hat{\underline{\chi}} - \frac{1}{2} \slashd \, \tr \underline{\chi} =0.
\end{eqnarray}
Hence $\tr \chi$ is a constant function over $S$ and we can assume $\tr \chi \equiv 1$ since we can always achieve this by modifying the null frame by a positive constant. So we get that
\begin{equation}
\slashd K=-\frac{1}{2}\slashdiv \hat{\underline{\chi}}.  \label{integrability condition}
\end{equation}
We consider the following operator $\slashdiv$ taking a 2-covariant symmetric, traceless tensor $\xi$ into the 1-form $\slashdiv \, \xi$. The $L^2$-adjoint of $\slashdiv$ is the operator taking a 1-form $f$ into the 2-covariant symmetric, traceless tensor $-\frac{1}{2}\widehat{\mathcal{L}_{f^{\sharp}} \eta|_{S}}$, where $\widehat{\mathcal{L}_{f^{\sharp}} \eta|_{S}}$ is the traceless part of the Lie derivative of $\eta|_{S}$ with respect to the vector field $f^{\sharp}$. This can be shown as the following:
\begin{eqnarray}
\int_S \langle \slashdiv \, \xi, f \rangle \, \dvol_{\eta|_S} 
&=& \int_S \langle \xi, - \slashnabla f \rangle \, \dvol_{\eta|_S} 
\\
&=&  \int_S \langle \xi, -\frac{1}{2}\sym{\slashnabla f} \rangle \,\dvol_{\eta|_S} 
\\
&=&  \int_{S} \langle \xi, -\frac{1}{2}\mathcal{L}_{f^{\sharp}} \eta|_S \rangle \, \dvol_{\eta_S}
\\
&=& \int_{S} \langle \xi, -\frac{1}{2}  \widehat{\mathcal{L}_{f^{\sharp}} \eta|_S} \rangle \,\dvol_{\eta_S}.
\end{eqnarray}
The kernel of the $L^2$-adjoint of $\slashdiv$ consists of the 1-form $f$ such that $f^{\sharp}$ is a conformal Killing vector field. Since the range of $\slashdiv$ is $L^2$ orthogonal to the kernel of its $L^2$-adjoint, then the identity \eqref{The generalized K-W identity} follows from \eqref{integrability condition}.


\begin{thebibliography}{100}

\bibitem{C} D. Christodoulou,
\emph{The Formation of Black Holes in General Relativity}.
EMS Monographs in Mathematics. European Mathematical Society (EMS), Z\"urich, 2009.

\bibitem{K-W} J. L. Kazdan, F. W. Warner,
Curvature Functions for Compact 2-Manifolds,
\emph{Ann. of Math.}  \textbf{99} (1974), 14--47.

\end{thebibliography}
\end{document}